\documentclass[11pt,titlepage]{amsart}
\usepackage{amsopn}
\usepackage{amssymb}
\usepackage{amsaddr}
\usepackage{graphicx}
\usepackage[all]{xy}
\usepackage{amscd}
\usepackage{color}
\usepackage[colorlinks]{hyperref}
\usepackage{amsmath,amsthm,amsfonts,amssymb}
\usepackage{verbatim}
\usepackage{latexsym}

 \newtheorem{theorem}{Theorem}[section]
\newtheorem{lemma}[theorem]{Lemma}
\newtheorem{proposition}[theorem]{Proposition}
\newtheorem{definition}[theorem]{Definition}

 \theoremstyle{definition}
 \theoremstyle{remark}

\begin{document}
\title[ R-duals of type $III$]
{ On R-duals of type $III$ in Hilbert spaces}

\author{Hartmut F\"{u}hr, Jahangir Cheshmavar and Ali Akbarnia}

\address{Hartmut F\"uhr, Lehrstuhl A f\"ur Mathematik, RWTH Aachen
University, D-52056 Aachen, Germany \\
\emph{fuehr@matha.rwth-aachen.de}}

\address{Jahangir Cheshmavar, Department of Mathematics,
Payame Noor University, P.O.Box 19395-3697, Tehran, Iran \\
\emph{j$_{_-}$cheshmavar@pnu.ac.ir}}

\address{Ali Akbarnia, Department of Mathematics,
Payame Noor University, P.O.Box 19395-3697, Tehran, Iran \\
\emph{Aliakbarnia7@gmail.com}}

\thanks{\it 2010 Mathematics Subject Classification: Primary
42C15; Secondary 53C05}

\keywords{Frames, Riesz sequence, Riesz basis, Spectral
representation, R-dual of type $I$, R-dual of type
$III$.}\dedicatory{} \commby{}
\begin{abstract}
Following work by Casazza, Kutyniok and Lammers, and extensions by
Stoeva and Christensen, we provide some
novel characterizations of $R$-dual sequences of type III in Hilbert spaces.
We systematically extend the construction procedure by basing it on a
choice of an antiunitary involution. For certain
classes of $R$-duals of type $III$, we derive a representation
of the associated frame operator in terms of spectral measures.

\end{abstract}
\maketitle


\section{Introduction}\label{Sec1}

In this paper we consider frame/Riesz sequence properties for a
sequence $\{f_i\}_{i\in I}$ in a Hilbert space $\mathcal{H}$ and
the corresponding sequence depending on two orthonormal bases, the
so-called  \textit{R-dual sequence (Riesz-dual sequence)}
$\{\omega_j\}_{j\in I}$ generated by combined action of some
operators on one of this orthonormal bases.
Sequences of this form were introduced in \cite{Casazza.2004} by
Casazza, Kutyniok and Lammers with the purpose of deriving duality principles for frames in arbitrary separable Hilbert space
$\mathcal{H}$. For each sequence in $\mathcal{H}$, they construct
a corresponding sequence depending on the choice of two orthonormal bases, with a
kind of duality relation between them, and then use them to derive duality principles for frames.
They called this sequence $R$-dual sequence. $R$-dual sequences have since been
considered by several authors, see
\cite{Christensen.2011,Christensen.2013,Chuang.2015}. One of the
most important extensions for the current paper appears in the work of Stoeva and
Christensen, which introduces various classes of $R$-duals, for example with the motivation to obtain
general versions of the duality principle in Gabor frames \cite{Christensen.2015}. In
\cite{Chuang.2015}, the authors give an equivalent condition of
two sequences to be $R$-duals.

In this paper we present variations of these constructions, that
are based on choices of certain isometric operators on $\ell^2(I)$ that help implement
the duality principles. Each choice of operator gives rise to a construction
of $R$-dual sequences that have the same desirable properties as the
original construction in \cite{Casazza.2004}, relating frame/Riesz
sequence properties of the sequence $\{f_i\}_{i\in I}$, to properties of its
$R$-dual $\{\omega_j\}_{j\in I}$, as studied in \cite{Casazza.2004},
but also by Christensen and Kim in \cite{Christensen.2011} and Stoeva
and Christensen in \cite{Christensen.2015}.
Finally, given a type $III$ of $R$-dual associated of a frame, we
derive a representation for its frame operator, (see Theorem
\ref{thm.10}) via the spectral theorem. The main results appear in Section 3 and Section 4.
Section 2 contains some basic definitions and results.


\section{Preliminaries}\label{Sec2}

In what follows, we will review basic definitions of frame and
Riesz basis and present type $I$ and $III$ of $R$-duals; for more
details, we refer the interested reader to the
\cite{Christensen.2008,Casazza.2001,Casazza.2004}. Throughout this
paper, let $\mathcal{H}$ be a separable Hilbert spaces,
$I_{\mathcal{H}}$ the identity operator on $\mathcal{H}$ and $I$
be a countable index set. Also let $F:=\{f_i\}_{i\in I}$ and
$\Omega:=\{\omega_j\}_{j\in I}$.\\

A collection of vectors $F$ in $\mathcal{H}$ is a \textit{Bessel
sequence} if there exists a constant $B>0$ so that
\begin{eqnarray*}
\sum_{i\in I}\mid\langle f, f_i\rangle \mid^2\leq B
\parallel f\parallel^2~,
\end{eqnarray*}
for all $f \in \mathcal{H}$. If in addition, there is a constant
$A>0$ so that
\begin{eqnarray*}
A \parallel f\parallel^2\leq \sum_{i\in I}\mid\langle f,
f_i\rangle \mid^2~,
\end{eqnarray*}
for all $f \in \mathcal{H}$, then $F$ is a \textit{frame} for
$\mathcal{H}$. The constants $A$ and $B$ are called frame bounds.
The frame $F$ is $A$-tight, if $A=B$. If at least the upper bounds
$B$ exists, $F$ is called a Bessel sequence. Here, the synthesis
operator of $F$ is defined by $$T_F: \ell^2(I)\rightarrow
\mathcal{H};\,\ T_F\{c_i\}_{i\in I}=\sum_{i\in I}c_if_i~.$$ Given
a frame $F$ in $\mathcal{H}$, its frame operator is
$$S_F:=T_FT_F^*:\mathcal{H}\rightarrow \mathcal{H};\,\
S_Ff=\sum_{i\in I}\langle f, f_i\rangle f_i~,$$ where  $T_F^*$ is
the adjoint of  $T_F$, and it is given by $T_F^*f=\{\langle f,
f_i\rangle\}_{i\in I}$. In this case $S_F$ is a bounded,
invertible, self-adjoint and positive operator. Moreover the
sequence $\widetilde{F}:=\{S^{-1}_Ff_i\}_{i\in I}$ is also a frame
for $\mathcal{H}$ satisfying the reconstruction formula
$f=\sum_{i\in I}\langle f, S^{-1}_Ff_i\rangle f_i$, for every $f
\in \mathcal{H}$. The sequence $\widetilde{F}$ is called the
canonical dual frame of $F$. Also, any sequence $G:=\{g_i\}_{i\in
I}$ in $\mathcal{H}$ which is not the canonical dual and satisfies
$f=\sum_{i\in I}\langle f, f_i\rangle g_i=\sum_{i\in I}\langle f,
g_i\rangle f_i$ is called an alternate dual frame of $F$.\\

A collection of vectors $\Omega$ in $\mathcal{H}$ is a
\textit{Riesz sequence} if there exist constants $C, D>0$ such
that
\begin{eqnarray*}
C\sum_{j\in I}\mid c_j\mid^2 \,\ \leq \,\ \parallel\sum_{j\in
I}c_j\omega_j\parallel^2\,\ \leq \,\ D\sum_{j\in I}\mid
c_j\mid^2~,
\end{eqnarray*}
for all finite sequences $\{c_j\}_{j\in I}$. The numbers $C$ and
$D$ are called Riesz bounds. A Riesz sequence $\Omega$ is a
\textit{Riesz basis} for $\mathcal{H}$ if
$\overline{span}\Omega=\mathcal{H}$.\\

We are now ready to introduce the main definitions considered in this
paper. We begin with the following well-known properties:
\begin{itemize}
  \item $F$ is a Bessel sequence in $\mathcal{H}$ if and only if $T_F^*$ is a well-defined and
  bounded operator on $\mathcal{H}$.
  \item $F$ is a frame for $\mathcal{H}$ if and only if $T_FoT_F^*:\mathcal{H}\rightarrow
  \mathcal{H}$ is a bounded and invertible operator.
  \item $F$ is a Riesz sequence in $\mathcal{H}$ if and only if $T_F^*oT_F:\ell^2(I) \rightarrow \ell^2(I)$
   is a bounded and invertible operator.
\end{itemize}
\begin{definition}
\indent Let $\{e_i\}_{i\in I}$ and $\{h_i\}_{i\in I}$ be
orthonormal bases for $\mathcal{H}$. Let $\{f_i\}_{i\in I}\subset
\mathcal{H}$ be such that $\sum_{i\in I}\mid\langle f_i,
e_j\rangle\mid^2<\infty, \forall j\in I$. In \cite{Casazza.2004}
the R-dual of $\{f_i\}_{i\in I}$ with respect to the orthonormal
bases $\{e_i\}_{i\in I}$ and $\{h_i\}_{i\in I}$ is defined as the
sequence given by
\begin{eqnarray*}
\omega_j=\sum_{i\in I}\langle f_i,e_j\rangle h_i, \quad\ j \in I.
\end{eqnarray*}

\end{definition}
This R-dual is called of type $I$ in \cite{Christensen.2015}. In
what follows, we will review type $III$ of R-duals, which are
essential to our main result in the next section.

\begin{definition}\cite{Christensen.2015}
Let $F$ be a frame for $\mathcal{H}$ with frame operator $S_F$.
Let $\{e_i\}_{i\in I}$ and $\{h_i\}_{i\in I}$ denote orthonormal
bases for $\mathcal{H}$ and $Q:\mathcal{H}\rightarrow\mathcal{H}$
be a bounded bijective operator with $\parallel
Q\parallel\leq\sqrt{\parallel S_F\parallel}$ and $\parallel
Q^{-1}\parallel\leq\sqrt{\parallel S^{-1}_F\parallel}$. The R-dual
of type $III$ of $F$ with respect to the triple $(\{e_i\}_{i\in
I}, \{h_i\}_{i\in I}, Q)$, is the sequence $\Omega$ defined by

\begin{eqnarray}
\omega_j=\sum_{i\in I}\langle S^{-1/2}_Ff_i,e_j\rangle Qh_i,
\quad\ j \in I.
\end{eqnarray}
\end{definition}

In this case $F$ obtained as

\begin{eqnarray} \label{10}
f_i=\sum_{j\in I}\langle \omega_j, (Q^*)^{-1}h_i\rangle
S^{1/2}_Fe_j, \quad\ \forall i \in I.
\end{eqnarray}

Relation (\ref{10}) does not imply, in general, that $F$ is an
R-dual of type $III$ of $\Omega$, that is, this definition of
R-dual is not symmetric. With appropriate choice of $Q$ the
symmetry property of the sequences $F$ and $\Omega$ characterized
in \cite[Theorems 4.4]{Christensen.2015}:\\

\begin{theorem}
\label{for2} Let $F$ be a frame for $\mathcal{H}$ and let $\Omega$
be a Riesz sequence with the same optimal bounds as $F$. Denote
the synthesis operator for $F$ by $T_F$ and the frame operators
for $F$ and $\Omega$ by $S_F$ and $S_{\Omega}$, respectively. If
$dim(ker T_F)=dim(span\{\omega_j\}^{\bot}_{j \in I})$, then there
exists orthonormal bases $\{e_i\}_{i\in I}$ and $\{h_i\}_{i\in I}$
for $\mathcal{H}$ such that
\begin{eqnarray}
\label{for3} \omega_j=\sum_{i\in I}\langle
f_i,S^{-1/2}_Fe_j\rangle \widetilde{S^{1/2}_{\Omega}} h_i, \quad\
j \in I~,
\end{eqnarray}
where $\widetilde{S^{1/2}_{\Omega}}$ is an extension of
$S^{1/2}_{\Omega}$ to an operator on $\mathcal{H}$.
\end{theorem}
The sequence $\Omega$ defined in (\ref{for3}) is called the
symmetrical R-dual of type $III$ of $F$ with respect to triple
$(\{e_i\}_{i\in I}, \{h_i\}_{i\in I},
\widetilde{S^{1/2}_{\Omega}})$. In this case, for all $i\in I$,
\begin{eqnarray}
\label{for.999}
 f_i=\sum_{j\in I}\langle
\omega_j,\widetilde{S^{-1/2}_{\Omega}}h_i\rangle S_F^{1/2}e_j~,
\end{eqnarray}
that is, $F$ is the R-dual of type $III$ of $\Omega$ with respect
to triple $(\{h_i\}_{i\in I}, \{e_i\}_{i\in I}, S^{1/2}_F)$.

\section{Characterizing r-duality in $\mathcal{H}$}\label{Sec3}

We first consider symmetrical R-duals of type $III$. In
\cite[Theorems 4.3]{Christensen.2015}, it is proved that if $F$ is
a frame sequence and $\Omega$ is an R-dual of $F$ of type $III$,
then the following hold:

\begin{itemize}
\item[(i)] $F$ is a frame for $\mathcal{H}$ if and only if
$\Omega$ is a Riesz sequence; in the affirmative case the bounds
for $F$ are also bounds for $\Omega$,\\
\item[(ii)] $F$ is a Riesz sequence if and only if $\Omega$ is a
frame for $\mathcal{H}$; in the affirmative case the bounds
for $F$ are also bounds for $\Omega$,\\
\item[(iii)] $\Omega$ is a Riesz Basis if and only if $F$ is a
Riesz Basis.
\end{itemize}
The following section presents another viewpoint of the construction of $R$-duals of type $III$.  In order to develop
this viewpoint, we introduce some additional terminology. A map $J: \mathcal{H} \to \mathcal{H}$ is called {\rm antiunitary} of $\mathcal{H}$
if $J$ is bijective, conjugate-linear, and isometric. It is called an {\em involution} if $J^2 = I_{\mathcal{H}}$.

The following lemma characterizes antiunitary involutions:
\begin{lemma}
 Let $\mathcal{H}$ denote a Hilbert space, and $J_0 : \mathcal{H} \to \mathcal{H}$ antiunitary involution.

 Let $J : \mathcal{H} \to \mathcal{H}$ denote a bijection. Then the following are equivalent:
 \begin{enumerate}
  \item[(a)] For all $u,v \in \mathcal{H}$, $\langle Ju, v \rangle = \overline{\langle u, Jv \rangle}$ and $J^2  = I_{\mathcal{H}}$
  \item[(b)] $J$ is conjugate-linear and satisfies $J^2 = I_{\mathcal{H}}$.
  \item[(c)] There exists a unitary map $U$ such that $J = U J_0$, and in addition $J_0 U J_0 = U^*$.
 \end{enumerate}
\end{lemma}
\begin{proof}
Assume (a)  For conjugate linearity of $J$, we observe for all $\lambda \in \mathbb{C}$ and all $u,v \in \mathcal{H}$ that by $(a)$,
 \[
  \langle J \lambda u, v \rangle = \overline{\langle \lambda u, Jv \rangle} = \overline{\lambda} \langle J u, v \rangle~,
 \] which establishes $J \lambda u = \overline{\lambda} J u$. By analogous reasoning, we get $J(u+v) = Ju + Jv$, and $(b)$ is established.

 For $(b) \Rightarrow (c)$, write $U = J J_0$. Then $U$ is isometric, bijective, and linear, hence unitary. Furthermore $U J_0 = J J_0^2 = J$, and $J^2 = J_0^2 = I_{\mathcal{H}}$ implies
 \[
  U J_0 U J_0 = J J_0 J_0 J J_0 J_0 = I_{\mathcal{H}}~,
 \] and thus $U^*  =  J_0 U J_0$.

 $(c) \Rightarrow (a)$: If $J = U J_0$ as in $(c)$, then $J^2 = U J_0 U J_0 U U^* = I_{\mathcal{H}}$. Since $J$ is bijective, conjugate-linear and isometric, polarization yields
 \[
  \langle J u, J w \rangle = \overline{\langle u,w \rangle}
 \] for all $u,v$, and with $w = Jv $ and $J^2 = I_{\mathcal{H}}$ we obtain
 \[
  \langle J u, v \rangle = \overline{ \langle u, J v \rangle }~,
 \]  which is $(a)$.
\end{proof}

We will explain below that any choice of antiunitary involution
gives rise to a construction of $R$-dual sequences that have all
the desirable properties of the original construction. The
original construction in \cite{Casazza.2004} is based on the
following concrete choice of $J_0: \ell^2(I) \to \ell^2(I)$
\[ J_0((x_i)_{i \in I}) = (\overline{x}_i)_{i \in I} ~.\] Since each Hilbert space is isometrically isomorphic to some $\ell^2(I)$, this example also establishes
that antiunitary involutions exist on every Hilbert space $\mathcal{H}$.

The next lemma notes how $J_0$, as defined above, enters in the definition of $R$-duals:
\begin{lemma}
\label{101}
 With the assumptions of the theorem (\ref{for2}), for
every $i,j \in I$,\,\
\begin{eqnarray}
\label{for.106}
\omega_j=T_{\widetilde{H}}oJ_0oT_{\widetilde{F}}^*(e_j)~,
\end{eqnarray}
where, $\widetilde{H}:=\{\widetilde{S_{\Omega}^{1/2}}h_i\}_{i \in
I},\,\ \widetilde{F}:=\{S^{-1/2}_Ff_i\}_{i\in I}$ and
\begin{eqnarray}
f_i=T_{\widetilde{E}}oJ_0oT_{\widetilde{\Omega}}^*(h_i)~,
\end{eqnarray}
 where,
$\widetilde{E}:=\{S^{1/2}_Fe_j\}_{j\in I}$,\,\
$\widetilde{\Omega}:=\{\widetilde{S_{\Omega}^{-1/2}}\omega_j\}_{j
\in I}$~.
\end{lemma}
\begin{proof}
By a direct and simple calculation we get (\ref{for.106}), which
is an operator representation of the symmetrical type $III$ of
R-duals $\Omega$ and $F$, stated in (\ref{for3}), respectively
(\ref{for.999}).
\end{proof}

Following this example, we define the symmetrical $R$-dual associated to an antiunitary involution $J$, a frame $F$ and a Riesz sequence $\Omega$ by replacing $J_0$ with $J$, leading to
\begin{eqnarray}
\label{for.106J}
\omega_j=T_{\widetilde{H}}oJoT_{\widetilde{F}}^*(e_j)~.
\end{eqnarray} Our aim is to prove that, for any choice of $J$, this construction yields a notion of $R$-dual that has all the desired properties of the original construction.

We will use the following simple Lemma.
\begin{lemma}
\label{lem.2} Let $U:\ell^2(I)\rightarrow \ell^2(I)$ is a linear
operator and $\widetilde{U}=JoUoJ$, where $J$ is an antiunitary involution. Then $\widetilde{U}$ is linear
with $\|U\|=\|\widetilde{U}\|$ and $\widetilde{U}^*=JoU^*oJ$.
\end{lemma}

\begin{proof}
$\langle \widetilde{U}f, g\rangle=\langle JoUoJf,
g\rangle=\overline{\langle UoJf, Jg\rangle}= \overline{\langle Jf,
U^*oJg\rangle}=\langle f, JoU^*oJg\rangle$ ~.
\end{proof}

The next result reads as follows:

\begin{proposition}
\label{102} Let $F$ is a Bessel sequence in $\mathcal{H}$,
$\mathcal{H}_F=\overline{span}F$, $P_F:\mathcal{H} \rightarrow
\mathcal{H}_F$ is an orthogonal projection. Then the following
statements hold:
\begin{itemize}
\item[(i)]$\widetilde{F}$ is a Riesz sequence if and only if
$P_FoP^*oPoP_F:\mathcal{H}_F \rightarrow \mathcal{H}_F$ is an
invertible operator, where $P:=S_F^{-1/2}$.

\item[(ii)]$T_{\Omega}^*=JoT_E^*oT_{\widetilde{F}}oJoT_{\widetilde{H}}^*$,
and
$T_{F}^*=JoT_H^*oT_{\widetilde{\Omega}}oJoT_{\widetilde{E}}^*~,$
where, $E:=\{e_i\}_{i\in I}$ and $H:=\{h_i\}_{i\in I}$. In
particular, $\Omega$ is a Bessel sequence in $\mathcal{H}$ if $F$
is a Bessel sequence in $\mathcal{H}$ and vice versa.
\end{itemize}
\end{proposition}

\begin{proof}
We first note that $\widetilde{F}:=\{Pf_i\}_{i\in I}$ is a Bessel
sequence in $\mathcal{H}$ with $T_{\widetilde{F}}^*=T_F^*oP^*$,
because for $f\in \mathcal{H}$,
\begin{eqnarray*}
\sum_{i\in I} |\langle f, Pf_i \rangle|^2\leq B\|P^*f\|^2=
B\|P\|^2\|f\|^2~,
\end{eqnarray*}
and $T_{\widetilde{F}}^*f=\{\langle f, Pf_i\rangle\}_{i\in I}$. In
particular, $\widetilde{F}$ is a Bessel sequence, if
$T_{\widetilde{F}}^*$ is a bounded operator on $\mathcal{H}$. On
the other hand for the proof of (i), we have
$$T_{\widetilde {F}}^*oT_{\widetilde{F}} =
T_F^*oP^*oPoT_F=T_F^*oP_FoP^*oPoP_FoT_F~,$$ and
$T_F:\ell^2(I)\rightarrow H_F$ is invertible. Hence,
$T_F^*:H_F\rightarrow \ell^2(I)$ is invertible; then
$\widetilde{F}$ is a Riesz sequence if and only if
$T_{\widetilde{F}}^*oT_{\widetilde{F}}$ is an invertible operator
on $\ell^2(I)$, if and only if $P_FoP^*oPoP_F$ is an invertible
operator on $\mathcal{H}_F$, as desired.

(ii) We have
\begin{eqnarray*}
T_{\Omega}^*f & = &\{\langle f, \omega_j\rangle\}_{j\in I}
=\{\langle f, T_{\widetilde{H}}oJoT_{\widetilde{F}}^*(e_j)
\rangle\}_{j\in I}=
 \{\langle T_{\widetilde{H}}^*f, JoT_{\widetilde{F}}^*(e_j) \rangle\}_{j\in
 I}\\
 & = & \{\overline{\langle JoT_{\widetilde{H}}^*f,
 T_{\widetilde{F}}^*(e_j)\rangle}\}_{j\in I}=
 \{\overline{\langle T_{\widetilde{F}}oJoT_{\widetilde{H}}^*f,
 e_j\rangle}\}_{j\in
 I}=\overline{T_E^*o(T_{\widetilde{F}}
 oJoT_{\widetilde{H}}^*f)}\\
 & = &(JoT_E^*oT_{\widetilde{F}}oJoT_{\widetilde{H}}^*)f~.
 \end{eqnarray*}
 That is,
 $T_{\Omega}^*=JoT_E^*oT_{\widetilde{F}}oJoT_{\widetilde{H}}^*$.
The process is similar for $T_{F}^*$.
\end{proof}
The next theorem shows that any antiunitary involution can be employed to define a notion of $R$-dual of type III.
\begin{theorem}
\label{thm.9}

Let $F$ is a frame sequence and $\Omega$ is an R-dual of type
$III$ of $F$ associated to the antiunitary involution $J$. Then the following hold:
\begin{eqnarray}
\label{for.098}
 T_{\Omega}oT_{\Omega}^* &=&
(T_{\widetilde{H}}oJoT_{\widetilde{F}}^*oT_EoJ)
o(JoT_E^*oT_{\widetilde{F}}oJoT_{\widetilde{H}^*})\\
&=&T_{\widetilde{H}}oJoT_{\widetilde{F}}^*oT_{\widetilde{F}}oJoT_{\widetilde{H}}^*~,
\end{eqnarray}
and
\begin{eqnarray}
\label{for.099}
T_{\Omega}^*oT_{\Omega}=(JoT_E^*oT_{\widetilde{F}}oJoT_{\widetilde{H}}^*)
o(T_{\widetilde{H}}oJoT_{\widetilde{F}}^*oT_EoJ)~.
\end{eqnarray}

 Furthermore,
\begin{itemize}
\item[(i)]$\Omega$ is a Riesz sequence in $\mathcal{H}$ if and
only if $F$ is a frame for $\mathcal{H}$.

\item[(ii)]$F$ is a Riesz sequence in $\mathcal{H}$ if and only if
$\Omega$ is a frame for $\mathcal{H}$.
\end{itemize}
\end{theorem}

\begin{proof}
Straightforward computation, using that $J^2=I_{\mathcal{H}}$,
$T_EoT_E^*=I_{\mathcal{H}}$, as well as Lemma (\ref{lem.2}) and Proposition
(\ref{102}), we have the formulas (\ref{for.098}) and
(\ref{for.099}). In particular,
\begin{itemize}
\item[(i)] $\Omega$ is a Riesz sequence if and only if
$T_{\Omega}^*oT_{\Omega}$ is a bounded and invertible operator on
$\ell^2(I)$
 if and only if
$$T_Fo\underbrace{JoT_{\widetilde{H}}^*oT_{\widetilde{H}}oJ}_{:=U}oT_F^*~,$$
is invertible, (which is by the fact that $S^{-1/2}_F$ is bounded
and invertible). Let now
$U:=JoT_{\widetilde{H}}^*oT_{\widetilde{H}}oJ$. Then
$U:\ell^2(I)\rightarrow \ell^2(I)$ is a positive and invertible
operator. Finally, for the proof of proposition, it is enough to
show that, $F$ is a frame for $\mathcal{H}$ if and only if
$T_FoUoT_F^*:\mathcal{H} \rightarrow \mathcal{H}$ is invertible.

Let $V:=U^{1/2}$, then $V$ is a positive and invertible operator
on $\ell^2(I)$ and we have
$$T_FoUoT_F^*=(VoT_F^*)^*o(VoT_F^*)~.$$
In particular, $F$ is a frame for $\mathcal{H}$ if and only if
$T_FoT_F^*:\mathcal{H}\rightarrow \mathcal{H}$ is a bounded and
invertible operator, if and only if $T_F^*:\mathcal{H} \rightarrow
\ell^2(I)$ is an embedding, if and only if $VoT_F^*:\mathcal{H}
\rightarrow \ell^2(I)$ is an embedding, if and only if
$(VoT_F^*)^*o(VoT_F^*)=T_FoUoT_F^*$ is an invertible operator, as
desired.

\item[(ii)] $F$ is a Riesz sequence in $\mathcal{H}$ if and only
if $T_{F}^*oT_{F}$ is a bounded and invertible operator on
$\ell^2(I)$ if and only if
$$T_{\Omega}o\underbrace{JoT_{\widetilde{E}}^*oT_{\widetilde{E}}oJ}_{:=W}oT_{\Omega}^*~,$$
is invertible, (which is by the fact that
$\widetilde{S_{\Omega}^{1/2}}$ is bounded and invertible). Let
$W:=JoT_{\widetilde{E}}^*oT_{\widetilde{E}}oJ$, with a similar
proof to $(i)$, $\Omega$ is a frame for $\mathcal{H}$.
\end{itemize}
\end{proof}

The generalization using antiunitary involutions highlights the algebraic and geometric properties of the map $J_0$ that
lead to the desirable properties of the $R$-dual. Further benefits of this additional freedom of choice in the design of $R$-duals remain
to be explored. For example, it is conceivable that a clever choice of $J$ may allow for $R$-duals that fulfill additional symmetry properties.

\section{Representation of $S_{\Omega}$}\label{Sec4}

In \cite{Chuang.2015}, Chuang and Zhao characterized $R$-duals $\Omega$ of
of type $I$ of a given frame $F$ by conditions, that are formulated without
explicit reference to the construction procedure of such duals.
In principle, this result shows that the frame poerator $S_\Omega$ is computable from $S_F$,
although the proof of the theorem is not explicit. The remainder of this section is
devoted to a more explicit description of $S_\Omega$ in terms of $S_F$. \\

\begin{theorem}\cite{Chuang.2015}
Let $F$ be a frame for $\mathcal{H}$ and $\Omega$ is a Riesz
sequence in $\mathcal{H}$. Denote the synthesis operator for $F$
by $T_F$, the frame operator of $F$ by $S_F$, and the frame
operator of $\Omega$ by $S_{\Omega}$. Then $\Omega$ is an R-dual
of type $I$ of $F$ if and only if the following two conditions hold:
\begin{itemize}
\item[(i)] there exists an antiunitary operator $\Lambda:
\mathcal{H}\rightarrow \overline{span}\Omega$ so that
$S_{\Omega}=\Lambda S_F \Lambda^{-1}$, \item[(ii)]
$dim(kerT_F)=dim((span\{\omega_j\}_{j\in I})^{\bot})$.
 \end{itemize}
\end{theorem}

In the following theorems, we present a representation for the
operator $S_{\Omega}$ associated with type $I$ and $III$ of
$R$-duals. Let now the countable index set $I$ be the integer
numbers set $\mathbb{Z}$. We use some ideas of
\cite{Casazza.2005}.

\begin{theorem} \label{for1000}
Let $F$ be a frame for $\mathcal{H}$ and $\Omega$ is a $R$-dual of
type $I$ of $F$ with respect to orthonormal bases $\{e_i\}_{i\in
I}$ and $\{h_i\}_{i\in I}$. Denote the frame operator of $F$ and
$\Omega$ by $S_F$ and $S_{\Omega}$, respectively, so that
\begin{eqnarray} \label{for1200} \sum_{j\in I}\mid \langle
S_{\Omega}h_0, h_j\rangle \mid < \infty~,
\end{eqnarray}
Then there exist operators $\{\mathcal{V}_j\}_{j\in I}$ on
$\mathcal{H}$ so that
\begin{itemize}
\item [(i)] \label{for100}
 $\mathcal{V}_j(S_{\Omega}h_0)=S_{\Omega}h_j~.$
\item [(ii)] If for all $k\in I, \,\ \{\mathcal{V}_j(h_k)\}_{j \in
I}$ is a Bessel sequence with Bessel bound $B$, then there exist
bounded operators $\{\Lambda_k\}_{k\in \mathbb{Z}}$, which satisfy
$\sup_{k\in \mathbb{Z}}\|\Lambda_k\|<\infty$, so that $S_{\Omega}$
associated with $F$ has the representation
\begin{eqnarray}
\label{for.10000} S_{\Omega}f=\sum_{i\in \mathbb{Z}} \langle f_i,
f_0\rangle \Lambda_i(f)~,
\end{eqnarray}
with unconditional convergence in the operator norm.
\end{itemize}
\end{theorem}

\begin{proof}
Let $\mathcal{U}:\mathcal{H}\rightarrow \mathcal{H}$ be unitary
shift operator $\mathcal{U}h_i:=h_{i+1}$, and let
$\mathcal{V}:\mathcal{H}\rightarrow \mathcal{H}$ be the operator
defined by $\mathcal{V}=S_{\Omega}\mathcal{U}S_{\Omega}^{-1}$.
Then we defined $\{\mathcal{V}_j\}_{j\in I}$ on $\mathcal{H}$ by
\begin{eqnarray*}
\mathcal{V}_j:=\mathcal{V}^j=S_{\Omega}\mathcal{U}^jS_{\Omega}^{-1}~.
\end{eqnarray*}
Then (i) is obviously fulfilled, since
\begin{eqnarray*}
\mathcal{V}_j(S_{\Omega}h_0)=\mathcal{V}^j(S_{\Omega}h_0)=
S_{\Omega}\mathcal{U}^jh_0=S_{\Omega}h_j, \quad\ \forall j\in I~.
\end{eqnarray*}
To proof of (ii), all calculations are similar to proof of theorem
4.5 in \cite{Casazza.2005}.
\end{proof}

For type $III$, the above procedure to computing $\langle
S_{\Omega}h_0, h_i\rangle$, does not work, and we have change our
strategy. First, recall (e.g., \cite{Kreyszig.1987}) that a representation of $S_{\Omega}$ in
terms of simple operators (projections) is the so-called a
\textit{spectral representation} of the operator $S_{\Omega}$.

 Let $T:\mathcal{H}\rightarrow \mathcal{H}$ be a bounded self-adjoint
linear operator. For $\lambda \in \mathbb{R}$, let
$T_{\lambda}=T-\lambda I_{\mathcal{H}}$. Denote the positive
square root of $T_{\lambda}^2$ by $|T_{\lambda}|$, the operator
$T_{\lambda}^+=\frac{1}{2}(|T_{\lambda}|+T_{\lambda})$ which is
called the positive part of $T_{\lambda}$, and the operator
$T_{\lambda}^-=\frac{1}{2}(|T_{\lambda}|-T_{\lambda})$ which is
called the negative part of $T_{\lambda}$. In this case,
$T_{\lambda}=T_{\lambda}^+-T_{\lambda}^-$, and the spectral family
$\mathcal{E}$ of $T$ is defined by
$\mathcal{E}=\{E_{\lambda}\}_{\lambda\in \mathbb{R}}$, where
$E_{\lambda}$ is the projection of $\mathcal{H}$ onto the null
space $\mathcal{N}(T_{\lambda}^+)$ of $T_{\lambda}^+$. Let
$$
m=\inf_{\|x\|=1}\langle Tx, x\rangle\,\ \mbox{and}\,\
M=\sup_{\|x\|=1}\langle Tx, x\rangle
$$
We now have the following theorem, \cite[Theorem
(9.2-1)]{Kreyszig.1987}:

\begin{theorem} \label{for2000}
Let $T:\mathcal{H}\rightarrow \mathcal{H}$ be a bounded
self-adjoint linear operator on a complex Hilbert space
$\mathcal{H}$. Then $T$ has the spectral representation as
follows:
\begin{eqnarray} \label{for2001}
T=mE_m+\int_m^M\lambda dE_{\lambda},
\end{eqnarray}
where the integral is to be understood in the sense of uniform
operator convergence. Also, for all $x, y \in \mathcal{H}$,
\begin{eqnarray} \label{for2002}
\langle Tx, y\rangle=mW(m)+\int_m^M\lambda d W(\lambda),
\end{eqnarray}
where, $W(\lambda)=\langle E_{\lambda}x, y\rangle$, and the
integral is an ordinary Riemann-Stieltjes integral.
\end{theorem}
Now we are ready to present a representation for the operator
$S_{\Omega}$ associated with type $III$ of $R$-duals.

\begin{theorem} \label{thm.10}
Let $F$ be a frame for $\mathcal{H}$ with frame operator $S_F$,
let $\Omega$ be a symmetrical $R$-dual of type $III$ of $F$, with
respect to the triple $(\{e_i\}_{i\in I}, \{h_i\}_{i\in I},
S_{\Omega}^{1/2})$, so that
\begin{eqnarray} \label{for1201} \sum_{j\in I}\mid \langle
S_{\Omega}h_0, h_j\rangle \mid < \infty~.
\end{eqnarray}
Assume there exist operators $\{\mathcal{V}_j\}_{j\in I}$ on
$\mathcal{H}$ such that
\begin{itemize}
\item [(i)] \label{for101}
 $\mathcal{V}_j(S_{\Omega}h_0)=S_{\Omega}h_j$, for all $j \in
 I$.
\item [(ii)] There is a constant $B>0$ so that for all $k\in I$
the set $\{\mathcal{V}_j(h_k)\}_{j \in I}$ is a Bessel sequence,
with Bessel bound $B$.
\end{itemize}

 Then there exist bounded operators $\{\Lambda_i\}_{i\in
\mathbb{Z}}$, which satisfy $\sup_{i\in
\mathbb{Z}}\|\Lambda_i\|<\infty$, and constants $\mathcal{C}_i$,
so that $S_{\Omega}$ associated with $F$ has the representation
\begin{eqnarray}
S_{\Omega}f=\sum_{i\in I}\left(M\langle
S_F^{-1/2}f_i,S_F^{-1/2}f_0 \rangle-\mathcal{C}_i\right)
\Lambda_i(f)~,
\end{eqnarray}
with unconditional convergence in the operator norm.
\end{theorem}

\begin{proof}
Since $S_{\Omega}$ is a bounded self-adjoint and linear operator on
$\mathcal{H}$, by the equation (\ref{for2002}) we have
\begin{eqnarray*}
\langle S_{\Omega} h_0, h_i\rangle & = & m W(m)+\int_m^M \lambda d
W(\lambda)\\ \mbox{(Riemann-Stieltjes integral)}&=&
mW(m)+MW(M)-mW(m)-\int_m^M W(\lambda)d\lambda,
\end{eqnarray*}
where $W(\lambda)=\langle E_{\lambda} h_0, h_i\rangle$ and
$E_{\lambda}$ is the projection of $\mathcal{H}$ onto the null
space $\mathcal{N}(S_{\Omega,\lambda}^+)$ of
$S_{\Omega,\lambda}^+$, with
$S_{\Omega,\lambda}=S_{\Omega}-\lambda I$. Therefore,

\begin{eqnarray}
\langle S_{\Omega} h_0, h_i\rangle =MW(M)-\int_m^M
W(\lambda)d\lambda
\end{eqnarray}

On the other hand, $\Omega$ is a Riesz sequence in $\mathcal{H}$
and then it is a Riesz basis for $V:=\overline{span}\Omega$. Then
the sequence $\{S_{\Omega}^{-1/2}\omega_j\}_{j\in I}$ is an
orthonormal basis for $V$ \cite[Lemma (1.1)]{Christensen.2015};
(note that $\mathcal{N}(S_{\Omega,\lambda}^+) \subset V)$.
Consider the extension $\widetilde{S_{\Omega}^{-1/2}}$ of
$S_{\Omega}^{-1/2}$ to an operator on $\mathcal{H}$, as in
\cite[Lemma (1.3)]{Christensen.2015};
$\{\widetilde{S_{\Omega}^{-1/2}}\omega_j\}_{j\in I}$ is an
orthonormal basis for $V$, too. Therefore, the orthogonal
projection $E_M$ of $\mathcal{H}$ onto
$\mathcal{N}(S_{\Omega,\lambda}^+)$ is given by
\begin{eqnarray}
\label{888} E_M f=\sum_{j\in I}\langle f,
\widetilde{S_{\Omega}^{-1/2}}\omega_j \rangle
\widetilde{S_{\Omega}^{-1/2}}\omega_j,\,\ \forall f \in
\mathcal{H}~.
\end{eqnarray}
It is enough to prove that if we define $E_M$ by (\ref{888}),
then $E_Mf=f$ for $f \in \mathcal{N}(S_{\Omega,\lambda}^+)$ and
$E_Mf=0$ for $f \in (\mathcal{N}(S_{\Omega,\lambda}^+))^{\perp}$;
the first equation follows by the orthonormality of
$\{S_{\Omega}^{-1/2}\omega_j\}_{j\in I}$, and the second by the
fact that the range of $S_{\Omega}^{-1/2}$ equals
$\mathcal{N}(S_{\Omega,\lambda}^+)$ because $S_{\Omega}^{-1/2}$ is
bijection on $\mathcal{N}(S_{\Omega,\lambda}^+)$. Therefore
\begin{eqnarray}
\langle E_M h_0, h_i\rangle &=& \sum_{j\in I}\langle h_0,
\widetilde{S_{\Omega}^{-1/2}}\omega_j \rangle \langle \widetilde{S_{\Omega}^{-1/2}}\omega_j, h_i\rangle\\
\mbox{ (by equation (\ref{for.999}))} &=& \langle S_F^{-1/2}f_i,
S_F^{-1/2}f_0\rangle~.
\end{eqnarray}
We now have
\begin{eqnarray*}
S_{\Omega}f&=&\sum_{j\in I}\langle f,h_j\rangle S_{\Omega}(h_j)\\
&=&\sum_{j\in I}\langle f,h_j\rangle \mathcal{V}_j\left(S_{\Omega}(h_0)\right)\\
&=&\sum_{j\in I}\langle f,h_j\rangle \mathcal{V}_j\left(\sum_{i\in
I}\langle S_{\Omega}(h_0),h_i\rangle h_i\right)\\
&=&\sum_{i\in I}\langle S_{\Omega}h_0,h_i\rangle\sum_{j\in
I}\langle f,h_j\rangle \mathcal{V}_j(h_i)\\
&=&\sum_{i\in I} \left(M\langle E_M h_0,
h_i\rangle-\mathcal{C}_i\right) \sum_{j\in
I}\langle f,h_j\rangle \mathcal{V}_j(h_i)\\
&=& \sum_{i\in I}\left(M\langle S_F^{-1/2}f_i,
S_F^{-1/2}f_0\rangle-\mathcal{C}_i\right) \sum_{j\in I}\langle
f,h_j\rangle \mathcal{V}_j(h_i)~,
\end{eqnarray*}
where $\mathcal{C}_i:=\int_m^M W(\lambda)d\lambda$. Defining
$\Lambda_i$ by $\Lambda_i(f)= \sum_{j\in I}\langle f,h_j\rangle
\mathcal{V}_j(h_i)$, then we have
\begin{eqnarray}
S_{\Omega}f=\sum_{i\in I}\left(M\langle
S_F^{-1/2}f_i,S_F^{-1/2}f_0 \rangle-\mathcal{C}_i\right)
\Lambda_i(f)~,
\end{eqnarray}

as desired. For convergence, the procedure is almost the same as
Theorem 4.5 in \cite{Casazza.2005}, but we are reviewing it. The
sequence $\{\mathcal{V}_j(h_i)\}_{j\in I}$ is a Bessel sequence
with bound $B$, then for any finite set $J\subset I$,
\begin{eqnarray*}
\|\sum_{j\in J}\langle f,h_j\rangle \mathcal{V}_j(h_i)\|^2 \leq B
\sum_{j\in J}\mid\langle f,h_j\rangle\mid^2=B\|f\|^2~,
\end{eqnarray*}
that is, the series defining $\Lambda_i$ are unconditionally
convergent and $\|\Lambda_i\|\leq\sqrt{B}$, $\forall k \in
\mathbb{Z}$.

 On the other hand, for finite subsets $I_1$ and $I_2$ of $I$, from the above
calculations,
\begin{eqnarray*}
\lefteqn{ \|\sum_{i\in I_1}\left(M\langle S_F^{-1/2}f_i,
S_F^{-1/2}f_0\rangle-\mathcal{C}_i\right) \sum_{j\in I_2}\langle
f,h_j\rangle \mathcal{V}_j(h_i)\|} \\ &=& \|\sum_{i\in I_1}\langle
S_{\Omega}h_0,h_i\rangle\sum_{j\in I_2}\langle
f,h_j\rangle \mathcal{V}_j(h_i)\| \\
  &\leq&  \left\{\sum_{i\in
I_1}|\langle S_{\Omega}h_0,h_i\rangle |\right\}\|\sum_{j\in
I_2}\langle f,h_j\rangle \mathcal{V}_j(h_i)\| \\
 &\leq& \left\{\sum_{i\in I_1}|\langle S_{\Omega}h_0,h_i\rangle
|\right\}\sqrt{B}\|f\|~.
\end{eqnarray*}
By the convergence of the series (\ref{for1200}), it follows that
the series in the construction of $\Lambda_i$ are unconditionally
convergent. Finally, for a finite subset $J\subset I$, we have
\begin{eqnarray*}
\lefteqn{ \|S_{\Omega}-\sum_{i\in J}\left(M\langle S_F^{-1/2}f_i,
S_F^{-1/2}f_0\rangle-\mathcal{C}_i\right) \Lambda_i\|} \\
&=&
\sup_{\|f\|=1}\|S_{\Omega}(f)-\sum_{i\in J}\left(M\langle
S_F^{-1/2}f_i, S_F^{-1/2}f_0\rangle-\mathcal{C}_i\right)
\Lambda_i(f)\|
\\ &=& \sup_{\|f\|=1}\|\sum_{i\in J^{c}}\left(M\langle S_F^{-1/2}f_i,
S_F^{-1/2}f_0\rangle-\mathcal{C}_i\right) \Lambda_i(f)\| \\
&\leq& \left\{\sum_{i\in J^{c}} \left|\left(M\langle
S_F^{-1/2}f_i,
S_F^{-1/2}f_0\rangle-\mathcal{C}_i\right) \right| \right\}\sup_{\|f\|=1}\sup_{i\in J^{c}}\|\Lambda_i(f)\| \\
&\leq& \left\{\sum_{i\in J^{c}}\left| \left(M\langle
S_F^{-1/2}f_i, S_F^{-1/2}f_0\rangle-\mathcal{C}_i\right| \right)
\right\} \sqrt{B}~.
\end{eqnarray*}
Now using (\ref{for1200}) and that $|\langle S_{\Omega}h_0,
h_i\rangle|=\left|\left(M\langle S_F^{-1/2}f_i,
S_F^{-1/2}f_0\rangle-\mathcal{C}_i\right)\right|$ , it follows
that the operators converge to $S_{\Omega}$ unconditionally in the
operator norm.
\end{proof}


\bibliographystyle{plain}
\end{document}